\documentclass[12pt]{article}  
\usepackage{amsthm,amsmath,amsfonts,epsfig,amssymb,latexsym}
\usepackage{mathrsfs}
\usepackage[T1]{fontenc}
\usepackage{graphicx}
\usepackage{enumerate}
\usepackage{xy,xypic}
\usepackage{graphics}
\usepackage{footnote}
\usepackage{scalefnt}
\usepackage{tikz}
\usetikzlibrary{shapes,arrows,fit,calc,positioning}
\usetikzlibrary{matrix}
\xyoption{all}
\title{{\Large On the complete intersection conjecture of Murthy}
} 
\author{
 Satya Mandal\footnote{Partially supported by a General Research Grant (no 2301857) from U. of Kansas}
 \\ 
{\small University of Kansas, Lawrence, Kansas 66045;}
{\small {\it  mandal@ku.edu 
  }}\\
 }  
\begin{document}
\renewcommand{\baselinestretch}{1.255}
\setlength{\parskip}{1ex plus0.5ex}
\date{7 October 2015} 
\newtheorem{theorem}{Theorem}[section]
\newtheorem{proposition}[theorem]{Proposition}
\newtheorem{lemma}[theorem]{Lemma}
\newtheorem{corollary}[theorem]{Corollary}
\newtheorem{construction}[theorem]{Construction}
\newtheorem{notations}[theorem]{Notations}
\newtheorem{question}[theorem]{Question}
\newtheorem{example}[theorem]{Example}
\newtheorem{definition}[theorem]{Definition} 
\newtheorem{clarification}[theorem]{Clarification} 
\newtheorem{Myproof}[theorem]{Proof} 
\newtheorem{remarks}[theorem]{Remarks}
\newtheorem{conjecture}[theorem]{Conjecture}
\newtheorem{openProblem}[theorem]{Open Problem}

\newcommand{\iso}{\stackrel{\sim}{\longrightarrow}}
\newcommand{\sur}{\twoheadrightarrow}
\newcommand{\bD}{\begin{definition}}
\newcommand{\eD}{\end{definition}}
\newcommand{\bP}{\begin{proposition}}
\newcommand{\eP}{\end{proposition}}
\newcommand{\bL}{\begin{lemma}}
\newcommand{\eL}{\end{lemma}}
\newcommand{\bT}{\begin{theorem}}
\newcommand{\eT}{\end{theorem}}
\newcommand{\bC}{\begin{corollary}}
\newcommand{\eC}{\end{corollary}} 
\newcommand{\eop}{\hfill \rule{2mm}{2mm}}
\newcommand{\pf}{\noindent{\bf Proof.~}}
\newcommand{\PD}{\text{proj} \dim}
\newcommand{\lra}{\longrightarrow}
\newcommand{\hra}{\hookrightarrow}
\newcommand{\llra}{\longleftrightarrow}
\newcommand{\Lra}{\Longrightarrow}
\newcommand{\Llra}{\Longleftrightarrow}
\newcommand{\bE}{\begin{enumerate}}
\newcommand{\eE}{\end{enumerate}}
\newcommand{\pic}{The proof is complete.}
\def\spec#1{\mathrm{Spec}\left(#1\right)}
\def\m{\mathfrak {m}}
\def\Sch{\underline{\mathrm{Sch}}}
\def\Sets{\underline{\mathrm{Sets}}}

\def\CA{\mathcal {A}}
\def\CB{\mathcal {B}}
\def\CP{\mathcal {P}}
\def\CC{\mathcal {C}}
\def\CD{\mathcal {D}}
\def\CE{\mathcal {E}}
\def\CF{\mathcal {F}}
\def\CE{\mathcal {E}}
\def\CG{\mathcal {G}}
\def\CH{\mathcal {H}}
\def\CI{\mathcal {I}}
\def\CJ{\mathcal {J}}
\def\CK{\mathcal {K}}
\def\CL{\mathcal {L}}
\def\CM{\mathcal {M}}
\def\CN{\mathcal {N}}
\def\CO{\mathcal {O}}
\def\CP{\mathcal {P}}
\def\CQ{\mathcal {Q}}
\def\CR{\mathcal {R}}
\def\CS{\mathcal {S}}
\def\CT{\mathcal {T}}
\def\CU{\mathcal {U}}
\def\CV{\mathcal {V}}
\def\CW{\mathcal {W}}
\def\CX{\mathcal {X}}
\def\CY{\mathcal {Y}}
\def\CZ{\mathcal {Z}}

\newcommand{\smallcirc}[1]{\scalebox{#1}{$\circ$}}
\def\BA{\mathbb {A}}
\def\BB{\mathbb {B}}
\def\BC{\mathbb {C}}
\def\BD{\mathbb {D}}
\def\BE{\mathbb {E}}
\def\BF{\mathbb {F}}
\def\BG{\mathbb {G}}
\def\BH{\mathbb {H}}
\def\BI{\mathbb {I}}
\def\BJ{\mathbb {J}}
\def\BK{\mathbb {K}}
\def\BL{\mathbb {L}}
\def\BM{\mathbb {M}}
\def\BN{\mathbb {N}}
\def\BO{\mathbb {O}}
\def\BP{\mathbb {P}}
\def\BQ{\mathbb {Q}}
\def\BR{\mathbb {R}}
\def\BS{\mathbb {S}}
\def\BT{\mathbb {T}}
\def\BU{\mathbb {U}}
\def\BV{\mathbb {V}}
\def\BW{\mathbb {W}}
\def\BX{\mathbb {X}}
\def\BY{\mathbb {Y}}
\def\BZ{\mathbb {Z}}

\newcommand{\TCP}{\textcolor{purple}}
\newcommand{\TCM}{\textcolor{magenta}}
\newcommand{\TCR}{\textcolor{red}}
\newcommand{\TCB}{\textcolor{blue}}
\newcommand{\TCG}{\textcolor{green}}

\def\SA{\mathscr {A}}
\def\SB{\mathscr {B}}
\def\SC{\mathscr {C}}
\def\SD{\mathscr {D}}
\def\SE{\mathscr {E}}
\def\SF{\mathscr {F}}
\def\SG{\mathscr {G}}
\def\SH{\mathscr {H}}
\def\SI{\mathscr {I}}
\def\SJ{\mathscr {J}}
\def\SK{\mathscr {K}}
\def\SL{\mathscr {L}}
\def\SN{\mathscr {N}}
\def\SO{\mathscr {O}}
\def\SP{\mathscr {P}}
\def\SQ{\mathscr {Q}}
\def\SR{\mathscr {R}}
\def\SS{\mathscr {S}}
\def\ST{\mathscr {T}}
\def\SU{\mathscr {U}}
\def\SV{\mathscr {V}}
\def\SW{\mathscr {W}}
\def\SX{\mathscr {X}}
\def\SY{\mathscr {Y}}
\def\SZ{\mathscr {Z}}

\def\bfA{{\bf A}}
\def\bfB{{\bf B}} 
\def\bfC{{\bf C}} 
\def\bfD{{\bf D}} 
\def\bfE{{\bf E}} 
\def\bfF{{\bf F}} 
\def\bfG{{\bf G}} 
\def\bfH{{\bf H}} 
\def\bfI{{\bf I}} 
\def\bfJ{{\bf J}} 
\def\bfK{{\bf K}} 
\def\bfL{{\bf L}} 
\def\bfM{{\bf M}} 
\def\bfN{{\bf N}} 
\def\bfO{{\bf O}} 
\def\bfP{{\bf P}} 
\def\bfQ{{\bf Q}} 
\def\bfR{{\bf R}} 
\def\bfS{{\bf S}} 
\def\bfT{{\bf T}} 
\def\bfU{{\bf U}} 
\def\bfV{{\bf V}} 
\def\bfW{{\bf W}} 
\def\bfX{{\bf X}} 
\def\bfY{{\bf Y}} 
\def\bfZ{{\bf Z}} 

\maketitle

\noindent{\bf Abstract}: {\it Suppose $A=k[X_1, X_2, \ldots, X_n]$ is a polynomial ring over a field $k$ and $I$ is an ideal in $A$. Then M. P. Murthy conjectured that $\mu(I)=\mu(I/I^2)$, where $\mu$ denotes the minimal number of generators. Recently, Fasel \cite{F} settled this conjecture, affirmatively, when  $k$ is an infinite perfect field, with $1/2\in k$ 
{\rm (always)}. We are able to do the same, when $k$ is an infinite field. In fact, we prove similar results for ideals $I$ in a polynomial ring $A=R[X]$, that contains a monic polynomial and $R$ is essentially finite type smooth algebra over an infinite field $k$, or $R$  is a regular ring over a perfect field $k$.
}

\section{Introduction}
One of the fundamental problems in commutative algebra, over last forty years,
 has been the  following conjecture of M. P. Murthy (\cite{M}, \cite[pp 85]{M1}), on complete intersections in affine spaces, as follows.

\begin{conjecture}[Murthy \cite{M}, \cite{M1}]\label{MurthyConj} 
Suppose $A=k[X_1, X_2, \ldots, X_n]$ is a polynomial ring over a field $k$.
Then, for any ideal $I$ in $A$, $\mu(I)=\mu(I/I^2)$, {\rm where $\mu$ denotes the minimal number of generators}.
\end{conjecture}
The conjecture (\ref{MurthyConj})  is sometimes  referred to as  Murthy's complete intersection conjecture, because if $I/I^2$ is free then it means $I$ would be a complete intersection
ideal. 

Recently, Fasel \cite{F} settled this conjecture, affirmatively,
when  $k$ is an infinite perfect field, with $1/2\in k$. We are able to do the same, when $k$ is an infinite field, with $1/2\in k$. 
In fact, we prove a much stronger theorem   (\ref{introMain}), given below.
A companion to Murthy's conjecture would be the following open problem.

\begin{openProblem}\label{OpenMonic}
Suppose $A=R[X]$ is a polynomial ring over a noetherian commutative ring $R$.
Suppose $I$ is an ideal in $A$ that contains a monic polynomial. Is $\mu(I)=\mu(I/I^2)$?
\end{openProblem}
For such an ideal $I$, as in (\ref{OpenMonic}), when $\mu(I/I^2) \geq \dim (A/I)+2$, Mohan Kumar (\cite{Mk}) proved that  $I$ is image of a projective $A$-module of rank $\mu(I/I^2)$
and it was proved in \cite{M2} that $\mu(I)=\mu(I/I^2)$. For  an ideal $I$, as in (\ref{OpenMonic}), without any further conditions, it was proved in \cite{MR}, that 
$I$ is set theoretically generated by $\mu(I/I^2)$-elements.

We settle this open problem (\ref{OpenMonic}), affirmatively, 
when $R$ is 
a regular ring, as specified below (\ref{introMain}).

%
\bT\label{introMain}
Let $R$ be a regular ring containing an infinite  field 
$k$, with $1/2\in k$. 
Assume $R$ is smooth  and essentially finite over $k$ or $k$ is perfect.
 Suppose $A=R[X]$ is the polynomial ring and
$I$ is an ideal in $A$ that contains a monic polynomial.
Then, $\mu(I)=\mu(I/I^2)$. 
In fact, for $n\geq 2$, any set of $n$-generators of $I/I^2$ lifts to a set of generators of $I$.
\eT
It was also well known that the validity of Murthy's conjecture would 
have implications on  the renowned  epi-morphism problem
(see \ref{ABH})
of S. Abhayankar \cite{A}.
For an exposition of the same the readers are referred to \cite{DG}, which would be our main reference on this. 
We state the epi-morphism  problem from \cite[Question 2.1]{DG}, as follows.

\begin{openProblem}[S.  Abhayankar]\label{ABH}
Suppose $k$ is a field, \\
 with $char(k)=0$. Let 
$$
\varphi:k[X_1, X_2, \ldots, X_n]\lra k[Y_1, \ldots, Y_m] \quad
{\rm be ~an~epi-morphism~of}~k{\rm-algebras}
$$
and $I=\ker(\varphi)$. Then,
 $I$ is generated by $n-m$ variables. That means $I=(F_1, \ldots, F_m)$ for some $F_1, \ldots, F_{n-m}\in I$ and 
$$
k[X_1, X_2, \ldots, X_n]=k[F_1, \ldots, F_{n-m}, Y_1', \ldots, Y_m'].
$$
 
\end{openProblem} 
A weaker version of the epi-morphism problem would be the following conjecture (see \cite[Question 2.2]{DG}).
\begin{conjecture}[S.  Abhayankar]\label{weakABH}
Suppose $k$ is a field, with $char(k)=0$ and 
$$
\varphi:k[X_1, X_2, \ldots, X_n]\lra k[Y_1, \ldots, Y_m] \quad
{\rm is ~an~epi-morphism~of}~k{\rm-algebras}
$$
and $I=\ker(\varphi)$. Then,
$\mu(I)=n-m$. 
\end{conjecture} 

As was indicated in \cite{DG}, very limited progress has been made on the Problem \ref{ABH}. Note that Conjecture \ref{weakABH} is subsumed  by Murthy's 
Conjecture \ref{MurthyConj} and same is true regarding the progress (see \cite{DG}).
A much stronger theorem, than the Conjecture \ref{weakABH}, follows from (\ref{introMain}) for polynomial algebras over regular rings $R$,
as specified below (\ref{introepiThm}). 
\bT\label{introepiThm}
Let $R$ be a regular ring over an infinite  field $k$, with $1/2\in k$.
Assume $R$ is smooth  and essentially finite over $k$ or $k$ is perfect.
Suppose 
$$
\varphi:R[X_1, X_2, \ldots X_n]\lra R[Y_1, Y_2, \ldots Y_m] \quad {\rm is~an~epi-morphism}
$$
of polynomial $R$-algebras and $I=\ker(\varphi)$. If $n-m\geq \dim R+1$, then $\mu(I)=\mu(I/I^2)$.
In particular, if $R$ is local, then $I$ is a complete intersection ideal.
\eT

Note that the hypotheses in (\ref{introMain}) entails only finite data (in fact, one equation). 
So, when $k$ is perfect, using the theorem of Popescu \cite{P}, for the purposes of the proofs of (\ref{introMain}),
we would be able to assume that $R$ is a smooth affine algebra over an infinite  field (see  the arguments in 
the proof of \cite[Theorem 2.1]{Sw}). With this approach, results in this article (e.g. \ref{FsMain327}) 
would be an improvement upon the respective versions in \cite{F}, relaxing the hypotheses in \cite{F} that the rings are of essentially finite type over 
infinite perfect fields.
Other than that,
we also give an alternate description (\ref{2pi0sEquiv}) of the obstruction set $\pi_0\left(Q_{2n}\right)(A)$, defined in \cite{F}.

 While the proof of Murthy's conjecture  in \cite{F} is 
elegant, its simplicity is even more astonishing. 
Proofs are further simplified, in this article, by proving that,
for   ideals $I$ in a polynomial ring 
$A=R[X]$, that contains a monic polynomial, 
 any local orientation  is homotopically trivial (\ref{trivialMonic}).
While reworking some of the proofs in \cite{F}, I also tried to elaborate and improve parts of the expositions
(e.g. \cite[Theorem 1.0.5]{F}) and  avoided  any
repetitions   that would be unwarranted.

%
We comment on the organization of this article. In \S \ref{secObsSets}, we elaborate the definition of the obstruction set $\pi_0\left(Q_{2n}\right)(A)$, and the 
obstruction classes and set up some notations. To avoid the stronger hypotheses in \cite{F}, 
in \S \ref{secFALilsting}, we restate and rework some of the results in \cite{F}.
In this section (\S \ref{secFALilsting}), we also record 
a statement of the homotopy lifting property theorem (\ref{MandalHT}), due to this author (unpublished), that was used in \cite{F}
and in (\ref{FsMain}).
%
In \S   \ref{SecMonicLift} we prove  the main theorem (\ref{introMain}), in this article.
In the subsection \S \ref{consequence}, we summarize the main of the consequences of (\ref{introMain}),
including the solution to Murthy's conjecture (\ref{MurthyConj}) and the weaker epi-morphism conjecture (\ref{weakABH}).
In \S  \ref{secALTObstru}, we provide an alternate description of the obstruction set $\pi_0\left(Q_{2n}\right)(A)$.

\vspace{5mm}
\noindent{\bf Acknowledgement}:{\it I would like to express my gratitude to 
M. P. Murthy for redrawing my interest in this problem and for subsequent discussions. 
In fact, it is because of his queries, I was convinced that there would be some interest for this article to others. 
I also thank Sankar P. Dutta for discussions on the work of Popescu.
}

\section{The Obstruction pre-sheaf}\label{secObsSets}
First, we establish some notations that  will be useful throughout this article.
\begin{notations}\label{nota}{\rm 
Throughout, $k$ will denote a field (or ring),
with $1/2\in k$ and $A, R$ will denote a commutative noetherian ring over $k$.
For a commutative ring $A$ and a finitely generated $A$-module $M$, the minimal number of 
generators of $M$ will be denoted by $\mu(M)$.

We denote 
$$
 q_{2n+1}=\sum_{i=1}^nX_iY_i+Z^2,\qquad \tilde{q}_{2n+1}=\sum_{i=1}^nX_iY_i+Z(Z-1).
 $$
Denote 
\begin{equation}
 \label{Q2nSept18}
 Q_{2n}=\spec{\SA_{2n}}~{\rm where}~\SA_{2n}=\frac{k[X_1, \ldots, X_n, Y_1, \ldots, Y_n, Z]}{\left(\tilde{q}_{2n+1}\right)} 
 \end{equation}
and 
\begin{equation}
 \label{Q2nSept18}
 Q_{2n}'=\spec{\SB_{2n}}~{\rm where}~\SB_{2n}=\frac{k[X_1, \ldots, X_n, Y_1, \ldots, Y_n, Z]}{\left(q_{2n+1}-1\right)}. 
 \end{equation}
There are inverse isomorphisms $\alpha: \SA_{2n}\iso \SB_{2n}$ 
 $\beta: \SB_{2n} \iso \SA_{2n}$ given by 
 \begin{equation}
 \label{alphaBeta}
 \left\{
 \begin{array}{ll}
 \alpha(x_i)=\frac{x_i}{2} & 1\leq i\leq n\\
  \alpha(y_i)=\frac{y_i}{2} & 1\leq i\leq n\\
\alpha(z)= \frac{z+1}{2}& \\
 \end{array}
 \right.
 \quad 
 \left\{
 \begin{array}{ll}
\beta(x_i)=2x_i & 1\leq i\leq n\\
\beta(y_i)=2y_i & 1\leq i\leq n\\
\beta(z)= 2z-1& \\
 \end{array}
 \right.
 \end{equation} 
Therefore, $Q_{2n}\cong Q_{2n}'$.
For a quadratic form $q$, of rank $n$, over a field $k$, $O(A, q)\subseteq GL_n(A)$ would denote the orthogonal 
 $q$-subgroup and $EO(A, q)\subseteq O(A, q)$ would denote the elementary orthogonal subgroup.  
 The category of schemes over $\spec{k}$ will be denoted by $\Sch_k$. Also, $\Sets$ will denote the category of sets.

}
\end{notations} 
The 
homotopy pre-sheaf $\pi_0\left(Q_{2n}\right)$ of sets,
 proved to be a key tool in \cite{F}, which we elaborate next.
 Recall, a contra-variant functor $\CF:\Sch_k \to \Sets$ is also called a pre-sheaf.

\bD{\rm 
Given a pre-sheaf $\CF:\Sch_k \to \Sets$, and a scheme $X\in \Sch_k$, define $\pi_0(\CF)(X)$ by the pushout 
\begin{equation}\label{Defpi0CF}
\diagram 
\CF(X\times \BA^1) \ar[r]^{T=0}\ar[d]_{T=1} & \CF(X)\ar[d]\\
\CF(X) \ar[r] & \pi_0(\CF)(X)\\ 
\enddiagram 
\qquad {\rm in}~~\Sets
\end{equation}
 For  an affine scheme $X=\spec{A}$ and a pre-sheaf $\CF$, as above, we write $\CF(A):=\CF(\spec{A})$ and $\pi_0(\CF)(A):=\pi_0(\CF)(\spec{A})$.
 So, $\pi_0(\CF)(A)$ is given by the pushout diagram 
 \begin{equation}\label{Defpi0Affine}
 \diagram 
 \CF(A[T]) \ar[r]^{T=0}\ar[d]_{T=1} & \CF(A)\ar[d]\\
\CF(A) \ar[r] & \pi_0(\CF)(A)\\ 
 \enddiagram
 \qquad {\rm in}~~ \Sets
 \end{equation}

 Given a scheme $Y\in \Sch_k$, the association
  $X\mapsto {\CH}om(X, Y)$ is a pre-sheaf on $\Sch_k$. 
 This pre-sheaf is often identified with $Y$,
  itself. So, in some  literature one may write, 
 $Y$ for the pre-sheaf ${\CH}om(-, Y)$ and 
 $Y(X):={\CH}om(X, Y)$.
Most importantly for us, for  schemes $X, Y\in \Sch_k$, the pre-sheaves  $\pi_0(Y)(X)$  are defined as in 
 diagram \ref{Defpi0CF}, or   \ref{Defpi0Affine}. 
 For the purposes of this article, $\pi_0\left(Q_{2n}\right)(X)$ and $\pi_0\left(Q_{2n}'\right)(X)$ would be of our particular interest. 
 For $X=\spec{A}$, it follows immediately that, $Q_{2n}(A)$ and $Q_{2n}'(A)$ can be identified with the sets, as follows:
  $$
 Q_{2n}(A)=\left\{(f_1,\ldots, f_n; g_1, \ldots, g_n; s)\in A^{2n+1}: \sum_{i=1}^nf_ig_i+s(s-1)=0 \right\}
 $$
  $$
 Q_{2n}'(A)=\left\{(f_1,\ldots, f_n;  g_1, \ldots, g_n; s)\in A^{2n+1}: \sum_{i=1}^nf_ig_i+s^2-1=0 \right\}
 $$
 The homotopy pre-sheaves are given by the pushout diagrams in $ \Sets$:
 $$
 \diagram 
  Q_{2n}(A[T]) \ar[r]^{T=0}\ar[d]_{T=1} & Q_{2n}(A)\ar[d]\\
Q_{2n}(A) \ar[r] & \pi_0\left(Q_{2n}\right)(A)\\ 
\enddiagram 
 \quad {\rm and }\quad 
  \diagram 
  Q_{2n}'(A[T]) \ar[r]^{T=0}\ar[d]_{T=1} & Q_{2n}'(A)\ar[d]\\
Q_{2n}'(A) \ar[r] & \pi_0\left(Q_{2n}'\right)(A)\\ 
\enddiagram 
 $$
  The isomorphism $Q_{2n}\cong Q_{2n}'$,  induces a bijection $\pi_0\left(Q_{2n}\right)(A) \cong \pi_0\left(Q_{2n}'\right)(A)$.

 For any ring $A$ and ${\bf v}=(f_1, \ldots, f_n; g_1, \ldots, g_n; s)\in Q_{2n}(A)$, let $I({\bf v}):=
 (f_1, \ldots, f_n, s)A$ denote the ideal. Also, let $\omega_{\bf v}:A^n\to \frac{I({\bf v})}{I({\bf v})^2}$ denote the 
 surjective map defined by assigning the standard basis $e_i\mapsto f_i+I^2$. Sometimes, $\omega_{\bf v}$ may be 
 called a local orientation. 
}
\eD
Now, we define local orientations of an ideal and the obstruction classes in $\pi_0\left(Q_{2n}\right)$.
\bD{\rm
Suppose $A$ is a commutative  ring and $I$ is an ideal in $A$. For an integer $n\geq 1$, 
a surjective homomorphism  $\omega:A^n \sur I/I^2$ would be called 
a \TCP{local $n$-orientation} of $I$. Clearly, a local $n$-orientation is determined by any
set of elements $f_1, \ldots, f_n\in I$ such that 
$I =(f_1, \ldots, f_n)+I^2$.  Further, given such a set of generators $f_1, \ldots, f_n$ of $I/I^2$, by Nakayama's lemma,
there is an $s\in I$ such that $(1-s)I\subseteq (f_1, \ldots, f_n)A$, and hence $\sum_{i=1}^nf_ig_i+s(s-1)=0$ for some 
$g_1, \ldots, g_n\in A$. Note,
$$
(f_1, \ldots, f_n; g_1, \ldots, g_n; s)\in Q_{2n}(A).
$$
Write
$$
\zeta(I, \omega):=[(f_1, \ldots, f_n; g_1, \ldots, g_n, s)]\in \pi_0\left(Q_{2n}(A)\right)
$$
It was established in \cite[Theorem 2.0.7]{F}, that this association is well defined. We refer to $\zeta(I, \omega)$,
as an 
\TCP{obstruction class}.  Therefore, we have a commutative diagram
$$
\diagram
Q_{2n}(A)\ar[rd]^{\zeta}\ar[d]_{\eta} & \\
\CO(A, n) \ar[r]_{\zeta} &  \pi_0\left(Q_{2n}(A)\right)\\
\enddiagram
 {\rm where}~\CO(A, n)={\rm  the~set~of~all~}n{\rm -orientations} ~(I, \omega)
$$
and $\eta({\bf v})=\left(I({\bf v}), \omega_{{\bf v}}\right)$.
Note, we use the same notation $\zeta$ for two set theoretic maps. 
}
\eD
\section{Homotopy and the lifting property}\label{secFALilsting}
In this section, we restate and rework, under the relaxed hypotheses in this article,
 some of the results in \cite{F} on the homotopy lifting property of the local orientations, to point out the modifications needed. 
First, we quote and interpret the following theorem from 
\cite{S}.
\bT\label{stravrova}
Suppose $A$ is a regular ring containing a perfect field $k$. Let $G$ be a reductive group
 scheme over $k$ such that every semi-simple  normal subgroup 
of $G$ contains $\BG_m^2$. 
Then, 
$$
\forall ~\sigma(T)\in G(A[T]), ~~\sigma(0)=1\Lra \sigma(T)\in E(A).
$$
In particular, 
with $G=O(q_{2n+1})$, we have 
$$
\forall ~\sigma(T)\in O(A[T], q_{2n+1})\subseteq GL_n(A[T]), ~~\sigma(0)=1\Lra \sigma(T)\in EO(A[T], q_{2n+1}).
$$
\eT 
\bC\label{WOperf}
Suppose $A$ is a regular ring containing a  field $k$. Then,
$$
\forall ~\sigma(T)\in O(A[T], q_{2n+1}), ~~\sigma(0)=1\Lra \sigma(T)\in EO(A[T], q_{2n+1}).
$$
\eC
\pf We will be following the arguments in \cite[Theorem 2.1]{Sw}, to reduce the problem to the perfect field case.
Let $\BF$ be the prime field of $k$. By including the coefficients of $\sigma(T)$, there is a finitely generated 
subalgebra $A' \subseteq A$ over $\BF$ such that $\sigma(T)\in A'[T]$. Note, $\sigma(T)\in O(A'[T], q_{2n+1})$ 
and $\sigma(0)=1$. Note $\BF \to A$ is geometrically regular. By Popescu's theorem \cite[Corollary 1.2]{Sw}, we have the diagram 
$$
\diagram 
\BF\ar[drr] \ar[r] & A'\ar[r]\ar@{^(->}[dr] & B\ar[d]\\
&& A\\
\enddiagram
~{\rm of~homomorphisms}
$$
such that $B$ is smooth over $\BF$. Since $\BF$ is perfect, by Theorem \ref{stravrova}, the $image(\sigma(T))$
in $B[T]$ is in $EO(B[T], q_{2n+1})$. Therefore, $\sigma(T)\in EO(A[T], q_{2n+1})$.
\pic$\eop$ 
%
\begin{remarks}{\rm
To reduce the problem to the perfect field case in (\ref{WOperf}), using Popescu's theorem \cite[Corollary 1.2]{Sw} 
has the same flavor of the similar reduction in \cite[pp 507]{T}, due to Mohan Kumar, which is more elementary.
One can work out an alternate proof of (\ref{WOperf}), using this argument of Mohan Kumar.
}
\end{remarks} 
For the benefit of the readers, 
we would elaborate the proof of \cite[Theorem 1.0.5]{F}, while reworking the following without the perfectness condition in \cite{F}.
\bT\label{1p0p5F}
Let $A$ be a smooth   and essentially of finite type algebra, 
over an infinite  field $k$, with $1/2\in k$. 
Then, for $n\geq 2$, the natural map
$$
\varphi: \frac{Q_{2n}'(A)}{EO(A, q_{2n+1})}\lra \pi_0\left(Q_{2n}'\right)(A)
\quad {\rm is~a~bijection}.
$$ 
\eT
\pf Clearly, $\varphi$ is well defined and is surjective. 
Now suppose $H(T)\in Q_{2n}'(A[T])$ is a homotopy. We need to show that there is a matrix 
$\tau\in EO(A, q_{2n+1})$ such that $H(0)=H(1)\tau$. 
%
For $R=\SB_{2n}, A[T], A$, use the following generic notations, to denote the quadratic modules
$$
\left\{
\begin{array}{lll}
q:=q_{2n+1}:R^{2n+1} \to R &{\rm sending} & (u_1, \ldots, u_n, v_1, \ldots, v_n, s)\mapsto \sum_{i=1}^nu_iv_i+s^2\\
q_0:R \to R &{\rm sending} & s\mapsto s^2\\
\end{array}
\right.
$$ 
As usual, define 
$B_q(e, e')=\frac{q(e+e')-q(e)-q(e')}{2}$. With respect to the standard basis, the matrix of $B_q$ is given by 
$$
B_q:=\frac{1}{2}\left(
\begin{array}{ccc} 
0 & I_n & 0\\ I_n & 0 & 0 \\ 0 & 0 & 2\\
\end{array}
\right)
$$
So, the map 
$
R^{2n+1} \to (R^{2n+1})^* ~{\rm sends} ~{\bf v} \mapsto {\bf v}B_q$. 
These bilinear forms give the following  exact sequences (write $\SB:=\SB_{2n}$)
$$
\diagram
0\ar[r] & \CK \ar[r] & \SB^{2n+1} \ar[rr]^{\langle ({\bf x}, {\bf y}, z), -\rangle} & & \SB \ar[r] & 0\\ 
0\ar[r] & K \ar[r] & A[T]^{2n+1} \ar[rr]^{\langle H(T), -\rangle} & & A[T] \ar[r] & 0\\ 
0\ar[r] & K_0 \ar[r] & A^{2n+1} \ar[rr]_{\langle H(0), -\rangle} & & A \ar[r] & 0\\ 
\enddiagram
$$
$$
{\rm So,}~~\CK=\left(\SB({\bf x}, {\bf y}, z)\right)^{\perp},\quad K=\left(A[T]H(T)\right)^{\perp},\quad  K_0=\left(RH(0)\right)^{\perp}
$$ 
are orthogonal complements, which inherit the  quadratic structures. 
It follows $\overline{K}:=K\otimes \frac{A[T]}{(T)}=\left(RH(0)\right)^{\perp}\cong K_0$. 
Therefore, there is an isometry,
$\sigma_0:K_0\iso \overline{K}$, which extends to 
$\sigma_0\otimes A[T]: K_0\otimes A[T]\iso \overline{K}\otimes A[T]$ an isometry.
({\it For clarity, note that it follows from Lindel's theorem (\cite{L}) that there is an isomorphism $\overline{K}\otimes A[T]\iso K$, 
which need not be an isometry}). 
Now, since $K$ is a torsor (see \cite{F}) that is locally trivial, by \cite[Theorem 3.1]{W}, 
there is, in deed, an isometry $\overline{K}\otimes A[T]\iso K$. 
By this identification, we say $\sigma_0$ extends to an isometry $\sigma_0\otimes A[T]:K_0\otimes A[T] \iso K$. 
Also, note 
$$
(A[T]H(T), q_{|A[T]H(T)}) \cong (A[T], q_0) \cong (A[T]H(0), q_{|A[T]H(0)})
$$
Putting all these together, there is an isometry $\sigma(T)\in O\left(A[T], q\right)$ such that the diagram 
\begin{equation}\label{linearMaps}
\diagram
0\ar[r] & K_0\otimes A[T] \ar[r]\ar[d]_{\sigma_0\otimes 1} & A[T]^{2n+1}\ar[d]^{\sigma(T)} \ar[rr]^{\langle H(0), -\rangle} & & A[T] \ar@{=}[d]\ar[r] & 0\\ 
0\ar[r] & K \ar[r] & A[T]^{2n+1} \ar[rr]_{\langle H(T), -\rangle} & & A[T] \ar[r] & 0\\ 
\enddiagram
\end{equation}
commutes. By construction, ({\it alternately, by composing with} $\sigma(0)^{-1}\otimes Id$), we have $\sigma(0)=I_{2n+1}$.
Now, by Corollary \ref{WOperf}, it follows $\sigma(T)\in EO(A[T])$ and hence $\sigma(1) \in EO(A)$. Since
$H(T)\sigma(T)=H(0)$, the proof is complete. 
 $\eop$ 

Before we proceed,  we define the action of $EO\left(A, q_{2n+1}\right)$ on $Q_{2n}(A)$ and 
give another  definition, for the convenience of subsequent discussions.
\bD{\rm 
Fix a commutative ring $A$. As usual, $EO\left(A, q_{2n+1}\right)$ acts on $A^{2n+1}$, which restricts to an action on $Q_{2n}'(A)$. Using the correspondences 
$\alpha: Q_{2n}(A) \iso Q_{2n}'(A)$, $\beta: Q_{2n}'(A) \iso Q_{2n}(A)$, define an action on $Q_{2n}(A)$ as follows:
$$
\forall~{\bf v}\in Q_{2n}(A), M\in EO\left(A,q_{2n+1}\right) ~~{\rm define} ~~{\bf v}*M:=\beta\left(\alpha({\bf v})M \right)
$$
This action is not  given by the usual matrix multiplication.  
Five different classes of the generators of  $EO\left(q_{2n+1}\right)(A)$  and their action on $Q_{2n}(A)$ are  given in \cite{F}.
}
\eD

\bD
Let $A$ be a commutative ring over $k$. Let ${\bf v}\in Q_{2n}(A)$. We write 
 ${\bf v}:=(a_1, \ldots, a_n; b_1, \ldots, b_n, s)$, 
 For integers, $r\geq 1$ we say that $r$-lifting property holds for ${\bf v}$, if 
 $$
 I({\bf v})=(a_1+\mu_1s^r, \ldots, a_n+\mu_ns^r) 
 \qquad {\rm for~some~}\mu_i\in A.
 $$ 
 We say the lifting property holds for ${\bf v}$, if 
 $$
 I({\bf v})=(a_1+\mu_1, \ldots, a_n+\mu_n) 
 \qquad {\rm for~some~}\mu_i\in I({\bf v})^2.
 $$ 
\eD

Before we allude to the key result in \cite[Corollary 3.2.6]{F} (see (\ref{FsMain}),
 we record the following homotopy lifting theorem, due to this author (unpublished),
  that was used crucially in the proof. 
\bT[Mandal]\label{MandalHT}
Let $R$ be a regular ring containing a field $k$. 
Let\\
 $H(T):=(f_1(T), \ldots, f_n(T), g_1(T), \ldots, g_n(T), s)\in Q_{2n}(R[T])$, with $s\in R$.\\
Write $a_i=f_i(0), b_i=g_i(0)$. 
Write $I(T)=(f_1(T), \ldots, f_n(T), s)$.
Also assume $I(0)=(a_1, \ldots, a_n)$ . 
Then,
$$
I(T)=(F_1, \ldots, F_n) \quad \ni ~~~f_i-F_i\in s^2R[T]
$$
\eT
\pf See \cite[Lemma 3.1.2]{F}. $\eop$ 
\bT\label{FsMain}
Suppose $A$ is a
regular ring containing  a filed $k$, with $1/2\in k$. Let $n\geq 2$ be an integer. 
%
 Let 
${\bf v} \in Q_{2n}(A)$ and $M\in EO\left(A, q_{2n+1}\right)$.
Then, ${\bf v}$ has $2$-lifting property if and only if ${\bf v}*M$ has the $2$-lifting property.
\eT
\pf We outline the  proof in \cite{F}. 
It would  be enough to assume that $M$ is a generator of $EO\left(A, q_{2n+1}\right)$.
There would be five cases to deal with, one for each type of generators of $EO\left(A, q_{2n+1}\right)$, listed in \cite[pp 3-4]{F}.
One of them, that is of the case of generators of the type 4 (in the list \cite[pp 3-4]{F}),  is fairly involved.
 This case follows from 
Theorem \ref{MandalHT} (see \cite[Lemma 3.1.2]{F}).  $\eop$

The following summarizes   the final results on homotopy and lifting of generators (also see \cite[Theorem 3.2.7]{F}).
\bT\label{FsMain327}
Suppose $A$ is a
regular ring containing an infinite  field  $k$, with $1/2\in k$. 
Assume $A$ is smooth  and essentially finite over $k$ or $k$ is perfect.
Let $n\geq 2$ be an integer.
 Denote ${\bf 0}:=(0,  \ldots, 0; 0, \ldots, 0;0)\in Q_{2n}(A)$ and let ${\bf v}\in Q_{2n}(A)$.
Then, the following conditions are equivalent:
\bE
\item \label{zetazero}The obstruction $\zeta\left(I({\bf v}, \omega_{\bf v})\right)=[{\bf 0}]\in \pi_0\left(Q_{2n}\right)(A)$.
\item\label{twoLift} ${\bf v}$ has $2$-lifting property.
\item\label{onlyLift} ${\bf v}$ has the lifting property.
\item\label{aNLift} ${\bf v}$ has $r$-lifting property, $\forall ~r\geq 2$.
\eE
\eT
\pf It is clear, $(\ref{twoLift}) \Lra (\ref{onlyLift})$. To prove $(\ref{onlyLift}) \Lra (\ref{zetazero})$,
 suppose $I({\bf v})=(a_1+\mu_1, \ldots, a_n+\mu_n)$,
with $\mu_i\in I({\bf v})^2$. 
Write ${\bf v}'=(a_1+\mu_1, \ldots, a_n+\mu_n; 0, \ldots, 0; 0)\in Q_{2n}(A)$. By \cite[2.0.10]{F}, we have
$
\zeta\left(I({\bf v}, \omega_{\bf v})\right)=\zeta\left(I({\bf v}', \omega_{{\bf v}'})\right)=[v_0]\in \pi_0\left(Q_{2n}\right).
$
This establishes, $(\ref{onlyLift}) \Lra (\ref{zetazero})$.

Now we prove $ (\ref{zetazero}) \Lra (\ref{twoLift})$. Assume $\zeta\left(I({\bf v}, \omega_{\bf v})\right)=[{\bf 0}]$.
In case $A$ is essentially finite over $k$, it follows from Theorem \ref{1p0p5F} that ${\bf 0}={\bf v}*M$, for some $M
\in EO(A, q_{2n+1})$ and (\ref{twoLift}) follows from Theorem \ref{FsMain}. However, when $A$ is regular and contains  an infinite perfect 
field, we have to use Popescu's theorem. 
By definition, $\zeta\left(I({\bf v}, \omega_{\bf v})\right)=[{\bf 0}]$ implies that there is a chain homotopy from ${\bf v}$ to ${\bf 0}$.
This data can also be encapsulated in a finitely generated algebra $A'$ over $k$. 
As in the proof of (\ref{WOperf}) there is a diagram
$$
\diagram 
k\ar[drr] \ar[r] & A'\ar[r]^{\iota}\ar@{^(->}[dr] & B\ar[d]\\
&& A\\
\enddiagram
~{\rm of~homomorphisms}
$$
such that $B$ is smooth over $k$. The homotopy relations are carried over to $B$. Therefore, by replacing $A$ by $B$, 
we can assume that $A$ is smooth and  essentially finitely  over $k$. So,  Theorem \ref{1p0p5F}  applies and 
(\ref{twoLift}) follows as in the previous case.


So, it is established that $(\ref{zetazero}) \Llra (\ref{twoLift}) \Llra (\ref{onlyLift})$. It is clear that 
$(\ref{aNLift}) \Lra (\ref{twoLift})$. Now suppose, one of the first three conditions hold. 
Fix $r\geq 2$. Notice $I({\bf v})=(a_1, \ldots, a_n, s^r)A$. So, replacement of $s$ by $s^r$ leads to the 
same obstruction class in  $\in \pi_0\left(Q_{2n}\right)(A)$, which is $=[{\bf 0}]\in \pi_0\left(Q_{2n}\right)(A)$. 
Since $(\ref{zetazero}) \Llra (\ref{twoLift})$, it follows 
$I({\bf v})$ has $2r$-lifting property and hence the $r$-lifting property.
\pic $\eop$

\section{Monic polynomials and the lifting property}\label{SecMonicLift}
In this section, first  we give an application of (\ref{FsMain327}), for ideals containing monic polynomials. 
In fact, we prove that all local orientations of such ideals  have trivial obstruction class, as follows.
\bP\label{trivialMonic}
Suppose  $R=A[X]$ is a polynomial ring over a commutative ring $A$ and $I$ is an ideal that contains a monic polynomial. 
Suppose $\omega: R^n\sur I/I^2$ is a surjective homomorphism ({\it local orientation}). 
Then $\zeta(I, \omega)=[{\bf 0}] \in \pi_0\left(Q_{2n}\right)(R)$, where 
${\bf 0}:=(0, 0, \ldots, 0,0, \ldots, 0)\in Q_{2n}(R)$.
%
\eP
\pf 
Let  $f_1, \ldots, f_n\in I$ be a lift of $\omega$. Then,
$
I= (f_1, f_2, \ldots, f_n)+I^2$.
We can assume that $f_1$ is a monic polynomial, with even degree and coefficient of the top degree term is $1$. 
Now, consider the transformation \cite{M2}:
$$
\varphi:A[X, T^{\pm1}]\iso A[X, T^{\pm1}]  ~~~{\rm by}~~~
\left\{
\begin{array}{l}
 \varphi(X)=X-T+T^{-1}\\
\varphi(T)=T\\
\end{array}\right. 
$$
There is a commutative diagram
$$
\diagram
A[X]\ar[d]\ar@{=}[r] & A[X]\\
A[X, T^{\pm1}] \ar[r] & A[X, T^{\pm1}]\ar[u]_{T=1}\\
\enddiagram
$$
Then, 
$\varphi(f_1)=f_1(X- T+T^{-1})$ is doubly monic in $T$, meaning that its lowest and the highest degree terms have coefficients $1$.
Let $F_1(X, T)=T^{\deg{f_1(X)}}\varphi(f_1)\in A[X, T]$. Then, $F_1(X, 0)=1$. Also, 
for $i=2, \ldots, n$ write $F_i(X, T)=T^{\delta}\varphi(f_i)$, for some integer $\delta \gg 0$, such that $F_i(X, T)\in TA[X, T]$.
Therefore,  $F_i(X,0)=0$. 
Now, write 
$$
\SI'=\varphi(IA[X, T^{\pm1}]) \quad {\rm and} \quad \SI:=\SI' \cap A[X, T].
$$
Since $\frac{A[X, T]}{\SI}\iso \frac{A[X, T^{\pm 1}]}{\SI'}$, it follows 
$$
 \SI=(F_1(X, T), \ldots, F_n(X, T))+\SI^2. 
 $$
Therefore, by Nakayama's Lemma, there is a $S(X, T)\in \SI$, such that  
$$
 (1-S(X, T))\SI\subseteq (F_1(X, T), F_2(X, T), \ldots , F_n(X, T).
$$
and hence
$$
\sum F_i(X, T) G_i(X, T)+S(X, T)(S(X, T)-1)=0 
 $$
for some $G_1, \ldots, G_n\in A[X, T]$. 
Write 
 $$
 \psi(X, T)= (F_1(X, T), F_2(X, T), \ldots, F_n(X, T); G_1(X, T), \ldots, G_n(X, T); S(X, T))
 $$
 Then, $\psi(X, T) \in Q_{2n}(A[X, T])$ and $\SI_{|T=1}=I$. Further, 
 $$
 \psi(X, 1)= (f_1, \ldots, f_n; G_1(X, 1), \ldots, G_n(X, 1); S(X,1))
 $$
 and
 $$
\psi(X, 0)=(1, 0, \ldots, 0; G_1(X, 0), \ldots, G_n(X, 0), S(X, 0)).
 $$
By \cite[2.0.10]{F}, $\psi(X,0) \sim {\bf 0}\in Q_{2n}(R)$. Hence, $\psi(X, 1)\sim {\bf 0}\in Q_{2n}(R)$.
Therefore,
$$
\zeta(I, \omega)=[ \psi(X, 1)]=[ {\bf 0}] \in \pi_0\left(Q_{2n}(R)\right). 
$$
\pic $\eop$ 

The following is our main theorem in this article, which is an extension of the main theorem in \cite{M2} mentioned in the introduction.
\bT\label{MoniMuImodI2}
Let $R$ be a regular ring over an infinite  field $k$, with $1/2\in k$ and $A=R[X]$
is the polynomial ring. 
Assume $R$ is smooth  and essentially finite over $k$ or $k$ is perfect.
Suppose, $I$ is an ideal in $A$ that contains a monic polynomial.
 Then  $\mu(I)=\mu(I/I^2)$. In fact, if $\mu(I/I^2)\geq 2$, any local orientation $\omega: A^n\sur I/I^2$ lifts to a 
set of generators of $I$.
\eT
\pf If $\mu(I/I^2)=1$, then $I$ is an invertible ideal with a monic polynomial, hence it is free.
The rest follows immediately from Proposition \ref{trivialMonic} and Theorem \ref{FsMain327}. \pic $\eop$

\subsection{The consequences}\label{consequence}
In this subsection, we summarize some of the consequences of the monic polynomial Theorem \ref{MoniMuImodI2}.
First, the following is the statement  on the solution of the conjecture of M. P. Murthy. The theorem is due to Fasel \cite[Theorem 3.2.9]{F}, 
in the case
when $k$ is perfect. However, our proof is much direct.
\bT\label{MurthyConj}
Let $k$ be an infinite  field, with $1/2\in k$. Let $A=k[X_1, \ldots, X_n]$ be the polynomial ring. 
Then, for any ideal $I$ of $A$, $\mu(I)=\mu(I/I^2)$.
\eT
\pf   By a change of variables (see \cite[Theorem 6.1.5]{M1}), we can assume that $I$ contains a monic polynomial
in $X_n$. Now, the proof is complete by (\ref{MoniMuImodI2}).  $\eop$

Following would be a more general version of Theorem \ref{MurthyConj}.
\bT\label{MoniMuNvariable}
Let $R$ be a regular ring over an infinite  field $k$, with $1/2\in k$ and $A=R[X_1,X_2, \ldots, X_n]$
is the polynomial ring in $n$ variables. 
Assume $R$ is smooth  and essentially finite over $k$ or $k$ is perfect.
Suppose, $I$ is an ideal in $A$ with $height(I)\geq \dim R+1$.
 Then  $\mu(I)=\mu(I/I^2)$. In fact, if $\mu(I/I^2)\geq 2$, any local orientation $\omega: A^n\sur I/I^2$ lifts to a 
set of generators of $I$.
\eT
\pf Again, by a change of variables (see \cite[Theorem 6.1.5]{M1}), we can assume that $I$ contains a monic polynomial
in $X_n$. Now, the proof is complete by (\ref{MoniMuImodI2}).  $\eop$

The following encompasses the solution of the weaker version of the epi-morphism conjecture (\ref{weakABH}), in the case when $k$ is an infinite fields. 
\bT\label{epiThmReg}
Let $R$ be a regular ring over an infinite  field $k$, with $1/2\in k$.
Assume $R$ is smooth  and essentially finite over $k$ or $k$ is perfect.
Suppose 
$$
\varphi:R[X_1, X_2, \ldots X_n]\lra R[Y_1, Y_2, \ldots Y_m] \quad {\rm is~an~epi-morphism}
$$
of polynomial $R$-algebras and $I=\ker(\varphi)$. If $n-m\geq \dim R+1$, then $\mu(I)=\mu(I/I^2)$.
In particular, if $R$ is local, then $I$ is a complete intersection ideal.
\eT 
\pf Since, $height(I)=n-m\geq \dim R+1$, it follows from (\ref{MoniMuNvariable}) that $\mu(I)=\mu(I/I^2)$. Note $I/I^2$
is a projective $R[Y_1, Y_2, \ldots Y_m]$-module of rank $m-n$. If $R$ is local, then $I/I^2$ is free or rank $n-m$ (see \cite[Theorem 2.1]{Sw}).
 Hence, $\mu(I)=mu(I/I^2)=n-m$. So, $I$ is a complete intersection ideal. 
\pic $\eop$

The following is the statement on the solution, in the case when $k$ is an infinite fields, of the weaker version of the epi-morphism conjecture (\ref{weakABH}).
\bC\label{epiWeak}
Suppose $k$ is an infinite field and 
$$
\varphi:k[X_1, X_2, \ldots X_n]\lra k[Y_1, Y_2, \ldots Y_m] \quad {\rm is~an~epi-morphism}
$$
of polynomial rings.  Then, $I$ is generated by $n-m$ elements.
\eC 
\pf Follows immediately from (\ref{epiThmReg}). 
$\eop$


\section{Alternate obstructions}\label{secALTObstru}
In this section, we give an alternate description of the obstruction sheaf $\pi_0\left(Q_{2n}\right)(A)$, which appears more traditional. 
\bD\label{Alterpi02n}{\rm 
Suppose $A$ is a commutative ring and $n\geq 1$ is an integer.
Write
$$
\CQ_n(A):=\left\{(f_1, \ldots, f_n, s)\in A^{n+1}: \exists ~g_1, \ldots, g_n\in A~\ni~\sum_{i=1}^nf_ig_i+s(s-1)=0  \right\}
$$
Note $A\mapsto \CQ_n(A)$ is a pre-sheaf on the category of affine schemes. As in diagram \ref{Defpi0CF}, one can define 
$\pi_0\left(\CQ_n(A)\right)$, by the pushout: 
 $$
 \diagram 
  \CQ_{n}(A[T]) \ar[r]^{T=0}\ar[d]_{T=1} & \CQ_{n}(A)\ar[d]\\
\CQ_{n}(A) \ar[r] & \pi_0\left(\CQ_{n}\right)(A)\\ 
\enddiagram 
 \qquad {\rm in}~~ \Sets.
 $$
 There is a natural map of sheaves $Q_{2n}(A) \to \CQ_n(A)$ sending 
 $$
 (f_1, \ldots, f_n; g_1, \ldots, g_n; s) \mapsto (f_1, \ldots, f_n, s)
 $$
 This induces a surjective map $\Phi: \pi_0\left(Q_{2n}\right(A) \sur \pi_0\left(\CQ_{n}\right)(A)$.
 }
\eD
\bL\label{2pi0sEquiv}
The map $\Phi:\pi_0\left(Q_{2n}\right)(A) \sur \pi_0\left(\CQ_{n}\right)(A)$ is a bijection. 
\eL 
\pf We define the inverse map $\Psi: \pi_0(\CQ_P) \lra \pi_0\left(Q_{2n}(A) \right)$.\\
Let $(f_1,\ldots, f_n, s) \in \CQ_n(A)$.  Then, 
$\sum f_ig_i+s(s-1)=0$ for some $g_i\in A$.
 We define,
$$
\Psi\left([(f_1, \ldots, f_n, s)]\right)= [(f_1,\ldots, f_n; g_1,\ldots, g_n; s)] \in \pi_0\left(Q_{2n}\right)(A).
$$
We need to show that $\Psi$ is well defined. With $I=(f_1, \ldots, f_n, s)$, let $\omega: A^n\sur \frac{I}{I^2}$
be induced by $f_1, \ldots, f_n$. By \cite[Lemma 2.0.10]{F}, the obstruction 
$\zeta(I, \omega_I):= [(f_1,\ldots, f_n; g_1,\ldots, g_n; s)]\in \pi_0\left(Q_{2n}(A) \right)$
 is independent of
$g_1, \ldots, g_n$. Therefore, $\Psi$ is well defined.
It is now clear that
$$
\Phi \Psi=Id\qquad {\rm and} \qquad
\Psi \Phi=Id. 
$$
\pic $\eop$
\end{document}